\documentclass{article}

\usepackage{amssymb,amsmath,amsthm}

\def\Dbar{\leavevmode\lower.6ex\hbox to 0pt{\hskip-.23ex
    \accent"16\hss}D}

\newtheorem{theorem}{Theorem}
\newtheorem{corollary}[theorem]{Corollary}
\newtheorem{definition}[theorem]{Definition}
\newtheorem{example}[theorem]{Example}
\newtheorem{problem}[theorem]{Problem}
\newtheorem{remark}[theorem]{Remark}

\newenvironment{keywords}{\begin{center}
\begin{minipage}[c]{11cm} {\bf Keywords:}}
{\end{minipage}
\end{center}}
\newenvironment{msc}{\begin{center}
\begin{minipage}[c]{11cm}
{\bf 2000 Mathematics Subject Classification:}}
{\end{minipage}
\end{center}
\bigskip}

\begin{document}

\title{Fractional Noether's theorem\\ in the Riesz-Caputo sense\thanks{Accepted 
(25/Jan/2010) for publication in \emph{Applied Mathematics and Computation}.}}

\author{Gast\~{a}o S. F. Frederico\\
        \texttt{gfrederico@mat.ua.pt}\\
        Department of Science and Technology\\
        University of Cape Verde\\
        Praia, Santiago, Cape Verde
        \and
        Delfim F. M. Torres\thanks{Corresponding author.
        Partially supported by the \emph{Centre for Research
        on Optimization and Control} (CEOC)
        of the University of Aveiro, cofinanced
        by the European Community fund FEDER/POCI 2010.}\\
        \texttt{delfim@ua.pt}\\
        Department of Mathematics\\
        University of Aveiro\\
        3810-193 Aveiro, Portugal}

\date{}

\maketitle

%--------------------------------------

\begin{abstract}
We prove a Noether's theorem for fractional variational problems
with Riesz-Caputo derivatives. Both Lagrangian and
Hamiltonian formulations are obtained. Illustrative examples
in the fractional context of the calculus of variations
and optimal control are given.
\end{abstract}


\begin{msc}
49K05, 26A33.
\end{msc}

%--------------------------------------

\section{Introduction}

Variational symmetries are defined by parameter transformations
that keep a problem of the calculus of variations or
optimal control invariant
\cite{MR2323264,MR0641031,CD:delfimPortMath04}.
Their importance, as recognized by Noether in 1918,
is connected with the existence of conservation laws
that can be used to reduce the order of the Euler-Lagrange differential equations
\cite{MR2351636,CD:GouveiaTorres:2005,MR1901565}.
Noether's symmetry theorem is nowadays recognized
as one of the most beautiful results
of the calculus of variations and optimal control
\cite{CD:Djukic:1972,CD:JMS:Torres:2002a,CD:delfimCPAA04}.

In 1967 a direct method for the problems of the calculus of variations,
which allow to obtain absolute extremizers directly,
without using the Euler-Lagrange equations,
was introduced by George Leitmann \cite{Leitmann67,Leit,Leitmann01}.
Time as shown that Leitmann's method is a general and fruitful principle
that can be applied with success to
a myriad of different classes of problems
\cite{MR2100262,CarlsonLeitmann05a,CarlsonLeitmann05b,UE,MR1954118,MR2065731,MR2035262}.
Interestingly, it turns out that Leitmann's and Noether's
principles are closely connected
\cite{SilvaTorres06,withLeitmann}.

The fractional calculus is an area of current strong
research with many different and important applications
\cite{Kilbas,CD:MilRos:1993,CD:Podlubny:1999,CD:SaKiMa:1993}.
In the last years its importance in the calculus of variations
and optimal control has been perceived, and a fractional variational
theory began to be developed by several different authors
\cite{AGRA6,CD:BalAv:2004,Cresson,MR2349759,MR2421931,CD:Hilfer:2000,CD:Muslih:2005,CD:Riewe:1997}.
Most part of the results in this direction make use of fractional derivatives
in the sense of Riemann-Liouville
\cite{CD:BalAv:2004,MR2349759,CD:FredericoTorres:2007,CD:FredericoTorres:2006,CD:Muslih:2005},
FALVA \cite{MR2421931,MR2256417,MR2405377},
or Caputo \cite{AGRA6,Porto2009:Ric,Caputo}. In 2007
generalized Euler-Lagrange fractional equations
and transversality conditions were studied for
variational problems defined in
terms of Riesz fractional derivatives \cite{CD:Agrawal:2007}.
In this paper we develop further the theory by obtaining
a fractional version of Noether's symmetry theorem
for variational problems with Riesz-Caputo derivatives
(Theorems~\ref{theo:TNfRL} and \ref{thm:mainResult:FDA06}).
Both fractional problems of the calculus
of variations and optimal control are considered.

We finish this introduction comparing in some details the results here obtained
with the ones of references \cite{MR2349759,MR2421931,MR2256417}.

In \cite{MR2349759} a $\left(\alpha,\beta\right)$ fractional derivative
is considered, which involves a right fractional derivative of
order $\alpha$ and a left fractional derivative of
order $\beta$ combined using a complex number $\gamma$.
The $\left(\alpha,\beta\right)$ fractional derivative
is useful when one needs to deal with complex valued functions.
In our paper we consider real valued functions only.
Moreover, we consider left and right derivatives in the sense
of Caputo, while the $\left(\alpha,\beta\right)$
derivative in \cite{MR2349759} is defined via left and right
Riemann-Liouville derivatives.
The advantage of using the Riesz symmetrized
Caputo fractional derivative instead of the
$\left(\alpha,\beta\right)$ derivative in \cite{MR2349759}
is that Caputo derivatives allow us to use
the standard boundary conditions of the calculus of variations,
which explains why they are more popular in engineer and physics.

In our paper we consider problems of the calculus
of variations for functions with one independent variable.
Paper \cite{MR2421931} initiates a new area of fractional
variational calculus by proposing
a fractional variational theory involving multiple integrals.
Some important consequences of such theory
in mechanical problems involving dissipative systems
with infinitely many degrees of freedom are given in \cite{MR2421931},
but a formal theory for that is missing.
Generalization of our present results to
multiple fractional variational integrals
is an interesting and challenging open question.
The recent results proved in \cite{MyID:182}
may be useful to that objective.

The results of \cite{MR2256417} are for
fractional Riemann-Liouville cost integrals
that depend on a parameter $\alpha$ but not on
fractional-order derivatives of order $\alpha$
as we do here: the variational problems of \cite{MR2256417}
are defined for Lagrangians that depend on the classical
derivative, while here we deal with fractional derivatives.

%--------------------------------------

\section{Preliminaries on Fractional Calculus}
\label{sec:fdRL}

In this section we fix notations by collecting the
definitions of fractional derivatives
in the sense of Riemann-Liouville, Caputo, and Riesz
\cite{CD:Agrawal:2007,CD:MilRos:1993,CD:Podlubny:1999,CD:SaKiMa:1993}.

\begin{definition}[Riemann-Liouville fractional integrals]
Let $f$ be a continuous function in the interval $[a,b]$. For $t
\in [a,b]$, the left Riemann-Liouville fractional integral
$_aI_t^\alpha f(t)$ and the right Riemann-Liouville fractional
integral $_tI_b^\alpha f(t)$ of order $\alpha$, $\alpha >0$,
are defined by
\begin{gather}
_aI_t^\alpha f(t) =
\frac{1}{\Gamma(\alpha)}\int_a^t (t-\theta)^{\alpha-1}f(\theta)d\theta \, ,
\label{eq:IFRLE} \\
_tI_b^\alpha f(t) = \frac{1}{\Gamma(\alpha)}\int_t^b
(\theta-t)^{\alpha-1}f(\theta)d\theta \, , \label{eq:IFRLD}
\end{gather}
where $\Gamma$ is the Euler gamma function.
\end{definition}

\begin{definition}[Riesz fractional integral]
Let $f$ be a continuous function in the interval $[a,b]$. For $t
\in [a,b]$, the Riesz fractional integral $_a^RI_b^\alpha f(t)$
of order $\alpha$, $\alpha > 0$, is defined by
 \begin{equation}
_a^RI_b^\alpha f(t) = \frac{1}{2\Gamma(\alpha)}\int_a^b
|t-\theta|^{\alpha-1}f(\theta)d\theta \, .
\label{eq:IFRz}
\end{equation}
\end{definition}

\begin{remark}
From equations \eqref{eq:IFRLE}--\eqref{eq:IFRz} it follows that
\begin{equation}\label{eq:demi}
_a^RI_b^\alpha f(t)=\frac{1}{2}\left[_aI_t^\alpha
f(t)+_tI_b^\alpha f(t)\right]\,.
\end{equation}
\end{remark}

\begin{definition}[fractional derivative in the sense of
Riemann-Liouville] Let $f$ be a continuous function
in the interval $[a,b]$. For $t \in [a,b]$, the left
Riemann-Liouville fractional derivative $_aD_t^\alpha f(t)$ and
the right Riemann-Liouville fractional derivative $_tD_b^\alpha
f(t)$ of order $\alpha$ are defined by
\begin{gather}
_aD_t^\alpha f(t)={D^n}_aI_t^{n-\alpha} f(t)=
\frac{1}{\Gamma(n-\alpha)}\left(\frac{d}{dt}\right)^{n}
\int_a^t (t-\theta)^{n-\alpha-1}f(\theta)d\theta \, , \label{eq:DFRLE} \\
_tD_b^\alpha f(t) ={(-D)^n}_tI_b^{n-\alpha}f(t)=
\frac{1}{\Gamma(n-\alpha)}\left(-\frac{d}{dt}\right)^{n}
\int_t^b(\theta - t)^{n-\alpha-1}f(\theta)d\theta \, ,
\label{eq:DFRLD}
\end{gather}
where $n \in \mathbb{N}$ is such that
$n-1 \leq \alpha < n$,
and $D$ is the usual derivative.
\end{definition}

\begin{definition}[fractional derivative in the sense of
Caputo] Let $f$ be a continuous function in $[a,b]$.
For $t \in [a,b]$, the left Caputo fractional derivative
$_a^CD_t^\alpha f(t)$ and the right Caputo fractional derivative
$_t^CD_b^\alpha f(t)$ of order $\alpha$ are defined in the
following way:
\begin{gather}
_a^CD_t^\alpha f(t)= {_aI_t^{n-\alpha}}D^nf(t)=
\frac{1}{\Gamma(n-\alpha)} \int_a^t
(t-\theta)^{n-\alpha-1}\left(\frac{d}{d\theta}\right)^{n}
f(\theta)d\theta \, , \label{eq:DFCE} \\
_t^CD_b^\alpha f(t) ={_tI_b^{n-\alpha}}(-D)^nf(t)=
\frac{1}{\Gamma(n-\alpha)} \int_t^b(\theta -
t)^{n-\alpha-1}\left(-\frac{d}{d\theta}\right)^{n}f(\theta)d\theta
\, , \label{eq:DFCD1}
\end{gather}
where $n \in \mathbb{N}$ is such that $n-1 \leq \alpha < n$.
\end{definition}

\begin{definition}[fractional derivatives in the sense of Riesz and Riesz-Caputo]
Let $f$ be a continuous function in $[a,b]$. For $t
\in [a,b]$, the Riesz fractional derivative $_a^RD_b^\alpha f(t)$
and the Riesz-Caputo fractional derivative $_a^{RC}D_b^\alpha
f(t)$ of order $\alpha$ are defined by
\begin{gather} _a^RD_b^\alpha f(t)=
D^n{_a^RI_t^{n-\alpha}}f(t)=
\frac{1}{\Gamma(n-\alpha)}\left(\frac{d}{dt}\right)^{n}
\int_a^b |t-\theta|^{n-\alpha-1}f(\theta)d\theta \, , \label{eq:DFRz} \\
_a^{RC}D_b^\alpha f(t) ={_a^RI_t^{n-\alpha}}D^nf(t)=
\frac{1}{\Gamma(n-\alpha)}
\int_a^b|t-\theta|^{n-\alpha-1}\left(\frac{d}{d\theta}\right)^{n}f(\theta)d\theta
\, , \label{eq:DFCz}
\end{gather}
where $n \in \mathbb{N}$ is such that $n-1 \leq \alpha < n$.
\end{definition}

\begin{remark}
Using equations \eqref{eq:demi} and
\eqref{eq:DFRLE}--\eqref{eq:DFCz} it follows that
\begin{equation*}
_a^RD_b^\alpha f(t)=\frac{1}{2}\left[_aD_t^\alpha
f(t)+(-1)^n\,{_tD_b^\alpha} f(t)\right]
\end{equation*}
and
\begin{equation*}
_a^{RC}D_b^\alpha f(t)=\frac{1}{2}\left[_a^CD_t^\alpha
f(t)+(-1)^n\,{_t^CD_b^\alpha} f(t)\right]\,.
\end{equation*}
In the particular case $0<\alpha<1$, we have:
\begin{equation}\label{eq:demi3}
_a^RD_b^\alpha f(t)=\frac{1}{2}\left[_aD_t^\alpha
f(t)-\,{_tD_b^\alpha} f(t)\right]
\end{equation}
and
\begin{equation}\label{eq:demi4}
_a^{RC}D_b^\alpha f(t)=\frac{1}{2}\left[_a^CD_t^\alpha
f(t)-\,{_t^CD_b^\alpha} f(t)\right]\,.
\end{equation}
\end{remark}

\begin{remark}
If $\alpha=1$, equalities \eqref{eq:DFRLE}--\eqref{eq:DFCD1}
give the classical derivatives:
\begin{equation*}
 _aD_t^1 f(t) =\, _a^CD_t^1 f(t)=\frac{d}{dt}
f(t) \, , \quad _tD_b^1 f(t) =\,_t^CD_b^1 f(t)= -\frac{d}{dt} f(t)
\, .
\end{equation*}
Substituting these quantities into \eqref{eq:demi3} and
\eqref{eq:demi4}, we obtain that
\begin{equation*}
 _a^RD_b^1 f(t) =\, _a^{RC}D_b^1 f(t)=\frac{d}{dt}
f(t) \, .
\end{equation*}
\end{remark}

%--------------------------------------

\section{Main Results}
\label{sec-ELRL}

In 2007 a formulation of the
Euler-Lagrange equations was given for problems
of the calculus of variations with fractional derivatives
in the sense of Riesz-Caputo \cite{CD:Agrawal:2007}.
Here we prove a fractional version of Noether's theorem
valid along the Riesz-Caputo fractional
Euler-Lagrange extremals \cite{CD:Agrawal:2007}. For that we
introduce an appropriate fractional operator that allow us
to generalize the classical concept of conservation law.
Under the extended fractional notion of conservation law
we begin by proving in \S\ref{sub:sec:CM} a fractional
Noether theorem without changing the time
variable $t$, \textrm{i.e.}, without transformation
of the independent variable (Theorem~\ref{theo:tnadf1}).
In \S\ref{sub:sec:NT} we proceed with a time-reparameterization
to obtain the fractional Noether's theorem in its general form
(Theorem~\ref{theo:TNfRL}). Finally, in \S\ref{sub:sec:OC}
we consider more general fractional optimal control problems
in the sense of Riesz-Caputo, obtaining
the corresponding fractional Noether's theorem in Hamiltonian form
(Theorem~\ref{thm:mainResult:FDA06}).

%--------------------------------------

\subsection{On the Riesz-Caputo conservation of momentum}
\label{sub:sec:CM}

We begin by defining the fractional functional under consideration.

\begin{problem}[The fractional problem of the calculus
of variations in the sense of Riesz-Caputo] The fractional problem
of the calculus of variations in the sense of Riesz-Caputo
consists to find the stationary functions of the functional
\begin{equation}
\label{Pf} I[q(\cdot)] = \int_a^b
L\left(t,q(t),{_a^{RC}D_b^\alpha} q(t)\right) dt \, ,
\end{equation}
where $[a,b] \subset \mathbb{R}$, $a<b$,
$0 < \alpha< 1$, and the admissible functions
$q: t \mapsto q(t)$ and the Lagrangian
$L : (t,q,v_l) \mapsto L(t,q,v_l)$ are assumed to be
functions of class $C^2$:
\begin{gather*}
q(\cdot) \in C^2\left([a,b];\,\mathbb{R}^n \right)\text{;}\\
L(\cdot,\cdot,\cdot) \in  C^2\left([a,b]\times\mathbb{R}^n \times
\mathbb{R}^n;\,\mathbb{R}\right)\text{.}
\end{gather*}
\end{problem}
Along the work, we denote by $\partial_{i}L$ the partial
derivative of $L$ with respect to its $i$-th argument, $i = 1, 2, 3$.
\begin{remark}
When $\alpha = 1$ the functional \eqref{Pf} is reduced to
the classical functional of the calculus of variations:
\begin{equation}
\label{eq:CFCV}
I[q(\cdot)] = \int_a^b L\left(t,q(t),\dot{q}(t)\right) dt \, .
\end{equation}
\end{remark}

The next theorem summarizes the main result of
 \cite{CD:Agrawal:2007}.

\begin{theorem}[\cite{CD:Agrawal:2007}]
\label{Thm:FractELeq1} If $q(\cdot)$ is an extremizer
of \eqref{Pf}, then it satisfies the
following \emph{fractional Euler-Lagrange
equation in the sense of Riesz-Caputo}:
\begin{equation}
\label{eq:eldf}
\partial_{2} L\left(t,q(t),{_a^{RC}D_b^\alpha q(t)}\right)
- {_a^RD_b^\alpha}\partial_{3} L\left(t,q(t),{_a^{RC}D_b^\alpha
q(t)}\right)  = 0
\end{equation}
for all $t \in [a,b]$.
\end{theorem}

\begin{remark}
The functional \eqref{Pf} involves Riesz-Caputo fractional
derivatives only. However, both Riesz-Caputo and Riesz fractional
derivatives appear in the fractional Euler-Lagrange equation
\eqref{eq:eldf}.
\end{remark}

\begin{remark}
Let $\alpha =1$. Then the fractional Euler-Lagrange
equation in the sense of Riesz-Caputo \eqref{eq:eldf}
is reduced to the classical Euler-Lagrange equation:
\begin{equation*}
\partial_{2} L\left(t,q(t),\dot{q}(t)\right)
-\frac{d}{dt}\partial_3 L\left(t,q(t),\dot{q}(t)\right)=0\, .
\end{equation*}
\end{remark}

Theorem~\ref{Thm:FractELeq1} leads to the concept of
fractional extremal in the sense of Riesz-Caputo.

\begin{definition}[fractional extremal in the sense of
Riesz-Caputo] \label{def:seC}
A function $q(\cdot)$ that is a solution of \eqref{eq:eldf}
is said to be a \emph{fractional Riesz-Caputo extremal}
for functional \eqref{Pf}.
\end{definition}

In order to prove a fractional Noether's theorem we adopt a
technique used in \cite{CD:FredericoTorres:2007,CD:Jost:1998}.
For that, we begin by introducing the notion of variational
invariance and by formulating a necessary condition
of invariance without transformation
of the independent variable $t$.

\begin{definition}[invariance of \eqref{Pf} without
transformation of the independent variable]
\label{def:inv1:MR}
Functional \eqref{Pf} is said to be invariant under an
$\varepsilon$-parameter group of infinitesimal
transformations
$\bar{q}(t)= q(t) + \varepsilon\xi(t,q(t)) + o(\varepsilon)$ if
\begin{equation}
\label{eq:invdf} \int_{t_{a}}^{t_{b}}
L\left(t,q(t),{_a^{RC}D_b^\alpha q(t)}\right) dt =
\int_{t_{a}}^{t_{b}} L\left(t,\bar{q}(t),{_a^{RC}D_b^\alpha
\bar{q}(t)}\right) dt
\end{equation}
for any subinterval $[{t_{a}},{t_{b}}] \subseteq [a,b]$.
\end{definition}

The next theorem establishes a necessary condition of invariance.

\begin{theorem}[necessary condition of invariance]
If functional \eqref{Pf} is invariant in the sense of
Definition~\ref{def:inv1:MR}, then
\begin{equation}
\label{eq:cnsidf}
\partial_{2} L\left(t,q(t),{_a^{RC}D_b^\alpha q(t)}\right) \cdot \xi(t,q(t))
+ \partial_{3} L\left(t,q(t),{_a^{RC}D_b^\alpha q(t)}\right) \cdot
{_a^{RC}D_b^\alpha \xi(t,q(t))} = 0 \, .
\end{equation}
\end{theorem}

\begin{remark}
Let $\alpha =1$. From \eqref{eq:cnsidf} we
obtain the classical condition of invariance of the
calculus of variations without transformation of the independent
variable $t$ (\textrm{cf.}, \textrm{e.g.}, \cite{Logan:b}):
\begin{equation*}
\label{eq:cnsiRT}
\partial_{2} L\left(t,q,\dot{q}\right)\cdot\xi(t,q)
+\partial_{3} L\left(t,q,\dot{q}\right)\cdot\dot{\xi}(t,q) = 0 \,
.
\end{equation*}
\end{remark}

\begin{proof}
Having in mind that condition \eqref{eq:invdf} is valid for any
subinterval $[{t_{a}},{t_{b}}] \subseteq [a,b]$, we can get rid
off the integral signs in \eqref{eq:invdf}. Differentiating this
condition with respect to $\varepsilon$, then substituting
$\varepsilon=0$, and using the definitions and properties of the
fractional derivatives given in Section~\ref{sec:fdRL}, we arrive
to the intended conclusion:
\begin{equation*}
\begin{split}
0
&= \partial_{2} L\left(t,q(t),{_a^{RC}D_b^\alpha q(t)}\right) \cdot \xi(t,q) \\
&+ \partial_{3} L\left(t,q(t),{_a^{RC}D_b^\alpha
q(t)}\right) \cdot
\frac{d}{d\varepsilon}\left[\frac{1}{\Gamma(n-\alpha)}
\int_a^b|t-\theta|^{n-\alpha-1}\left(\frac{d}{d\theta}\right)^{n}\bar{q}(\theta)d\theta
\right]_{\varepsilon=0}\\
&= \partial_{2} L\left(t,q,{_a^{RC}D_b^\alpha q}\right)\cdot\xi(t,q) \\
&\qquad+ \partial_{3} L\left(t,q,{_a^{RC}D_b^\alpha q}\right)
\cdot\frac{d}{d\varepsilon} \left[\frac{1}{\Gamma(n-\alpha)}
\int_a^b|t-\theta|^{n-\alpha-1}\left(\frac{d}{d\theta}\right)^{n}q(\theta)d\theta
\right. \\
&\qquad+\left.\frac{\varepsilon}{\Gamma(n-\alpha)}\int_a^b
|t-\theta|^{n-\alpha-1}\left(\frac{d}{d\theta}\right)^{n}\xi(\theta,q)d\theta\right]_{\varepsilon=0}  \\
& = \partial_{2} L\left(t,q,{_a^{RC}D_b^\alpha q}\right)\cdot\xi(t,q) \\
&\qquad+ \partial_{3} L\left(t,q,{_a^{RC}D_b^\alpha q}\right)
\cdot\frac{1}{\Gamma(n-\alpha)}\int_a^b
|t-\theta|^{n-\alpha-1}\left(\frac{d}{d\theta}\right)^{n}\xi(\theta,q)d\theta\ \\
&=\partial_{2} L\left(t,q,{_a^{RC}D_b^\alpha q}\right) \cdot
\xi(t,q) + \partial_{3} L\left(t,q,{_a^{RC}D_b^\alpha q}\right)
\cdot {_a^{RC}D_b^\alpha} \,\xi(t,q) \,.
\end{split}
\end{equation*}
\end{proof}

The following definition is useful in order to introduce an
appropriate concept of \emph{fractional conserved quantity in the
sense of Riesz-Caputo}.

\begin{definition}
\label{def:oprl}
 Given two functions $f$ and $g$ of class
$C^1$ in the interval $[a,b]$, we introduce the following
operator:
\begin{equation*}
\mathcal{D}_{t}^{\gamma}\left(f,g\right) = g \cdot {_a^RD_b^\gamma} f
+ f \cdot {_a^{RC}D_b^\gamma} g \, ,
\end{equation*}
where $t \in [a,b]$ and $\gamma \in \mathbb{R}_0^+$.
\end{definition}

\begin{remark}
Similar operators were used in
\cite[Definition~19]{CD:FredericoTorres:2007}
but involving Riemann-Liouville fractional derivatives.
We note that the new operator $\mathcal{D}_{t}^{\gamma}$
proposed here involves both Riesz and Riesz-Caputo
fractional derivatives.
\end{remark}

\begin{remark}
\label{rem:oprl}
In the classical context one has $\gamma=1$ and
$$
\mathcal{D}_{t}^{1}\left(f,g\right)=f' \cdot g + f \cdot g' \\
= \frac{d}{dt}(f \cdot g)=\mathcal{D}_{t}^{1}\left(g,f\right) \, .
$$
Roughly speaking, $\mathcal{D}_{t}^{\gamma}\left(f,g\right)$
is a fractional version of the derivative of the product
of $f$ with $g$. Differently from the classical context,
in the fractional case one has, in general,
$\mathcal{D}_{t}^{\gamma}\left(f,g\right)
\ne \mathcal{D}_{t}^{\gamma}\left(g,f\right)$.
\end{remark}

We now prove the fractional Noether's theorem in the
sense of Riesz-Caputo without transformation
of the independent variable $t$.

\begin{theorem}[Noether's theorem in the
sense of Riesz-Caputo without transformation of time]
\label{theo:tnadf1} If functional \eqref{Pf} is invariant in the
sense of Definition~\ref{def:inv1:MR}, then
\begin{equation}
\label{eq:LC:Frac:RL}
\mathcal{D}_{t}^{\alpha}\left[\partial_{3}
L\left(t,q(t),{_a^{RC}D_b^{\alpha} q(t)}\right),\xi(t,q(t))\right] = 0
\end{equation}
along any fractional Riesz-Caputo extremal $q(t)$, $t \in [a,b]$
(Definition~\ref{def:seC}).
\end{theorem}

\begin{remark}
\label{rem:22}
In the particular case when $\alpha = 1$
we get from the fractional conservation law
in the sense of Riesz-Caputo
\eqref{eq:LC:Frac:RL} the classical
Noether's conservation law of momentum
(\textrm{cf.}, \textrm{e.g.}, \cite{CD:Jost:1998,Logan:b}):
\begin{equation*}
\frac{d}{dt} \left[ \partial_{3}
L\left(t,q(t),\dot{q}(t)\right)\cdot\xi(t,q(t)) \right] = 0
\end{equation*}
along any Euler-Lagrange extremal $q(\cdot)$ of \eqref{eq:CFCV}.
For this reason, we call the fractional law \eqref{eq:LC:Frac:RL}
\emph{the fractional Riesz-Caputo conservation of momentum}.
\end{remark}

\begin{proof}
Using the fractional Euler-Lagrange equation \eqref{eq:eldf},
we have:
\begin{equation}
\label{eq:eldf3}
\partial_{2} L\left(t,q,{_a^{RC}D_b^\alpha q}\right)
= {_{a}^{R}D_b^\alpha}\partial_{3} L\left(t,q,{_a^{RC}D_b^\alpha
q}\right) \, .
\end{equation}
Replacing \eqref{eq:eldf3} in the necessary
condition of invariance \eqref{eq:cnsidf}, we get:
\begin{equation}
\label{eq:dtnrl} {_{a}^{R}D_b^\alpha}\partial_{3}
L\left(t,q,{_a^{RC}D_b^\alpha q}\right)\cdot\xi(t,q) +\partial_{3}
L\left(t,q,{_a^{RC}D_b^\alpha q}\right)\cdot{_a^{RC}D_t^\alpha
\xi(t,q)}=0\,.
\end{equation}
By definition of the operator
$\mathcal{D}_{t}^{\gamma}\left(f,g\right)$ it results from
\eqref{eq:dtnrl} that
\begin{equation*}
\mathcal{D}_{t}^{\alpha}\left[\partial_{3}
L\left(t,q,{_a^{RC}D_b^\alpha q}\right),\xi(t,q)\right] = 0 \, .
\end{equation*}
\end{proof}

%--------------------------------------

\subsection{The Noether theorem in the sense of Riesz-Caputo}
\label{sub:sec:NT}

The next definition gives a more general notion
of invariance for the integral functional \eqref{Pf}.
The main result of this section, the
Theorem~\ref{theo:TNfRL}, is formulated
with the help of this definition.

\begin{definition}[invariance of \eqref{Pf}]
\label{def:invadf} The integral functional \eqref{Pf}
is said to be
invariant under the one-parameter group of infinitesimal
transformations
\begin{equation*}
\begin{cases}
\bar{t} = t + \varepsilon\tau(t,q(t)) + o(\varepsilon) \, ,\\
\bar{q}(t) = q(t) + \varepsilon\xi(t,q(t)) + o(\varepsilon) \, ,
\end{cases}
\end{equation*}
if
\begin{equation*}
\int_{t_{a}}^{t_{b}} L\left(t,q(t),{_a^{RC}D_b^\alpha q(t)}
\right) dt = \int_{\bar{t}(t_a)}^{\bar{t}(t_b)}
L\left(\bar{t},\bar{q}(\bar{t}), {_a^{RC}D_b^\alpha
\bar{q}(\bar{t})}\right) d\bar{t}
\end{equation*}
for any subinterval $[{t_{a}},{t_{b}}] \subseteq [a,b]$.
\end{definition}

Our next theorem gives a generalization
of Noether's theorem for fractional problems
of the calculus of variations
in the sense of Riesz-Caputo.

\begin{theorem}[Noether's fractional theorem in the sense of
Riesz-Caputo] \label{theo:TNfRL} If the integral functional
\eqref{Pf} is invariant in the sense of
Definition~\ref{def:invadf}, then
\begin{multline}
\label{eq:tndf}
\mathcal{D}_{t}^{\alpha}\left[\partial_{3} L\left(t,q,{_a^{RC}D_t^\alpha q}\right),
\xi(t,q)\right] \\
+ \mathcal{D}_{t}^{\alpha}\left[L\left(t,q,{_a^{RC}D_t^\alpha q}\right) -
\alpha\partial_{3} L\left(t,q,{_a^{RC}D_t^\alpha
q}\right)\cdot{_a^{RC}D_b^\alpha q}, \tau(t,q)\right] = 0
\end{multline}
along any fractional Riesz-Caputo extremal $q(\cdot)$.
\end{theorem}

\begin{remark}
\label{rem:25}
In the particular case $\alpha= 1$ we obtain from \eqref{eq:tndf}
the classical Noether's conservation law (\textrm{cf.}, \textrm{e.g.},
\cite{CD:Jost:1998,Logan:b}):
\begin{equation*}
\frac{d}{dt}
\left[\partial_{3} L\left(t,q,\dot{q}\right)\cdot\xi(t,q)
+ \left( L(t,q,\dot{q}) - \partial_{3} L\left(t,q,\dot{q}\right)
\cdot \dot{q} \right) \tau(t,q)\right] = 0
\end{equation*}
along any Euler-Lagrange extremal $q(\cdot)$ of \eqref{eq:CFCV}.
\end{remark}

\begin{proof}
Our proof is an extension of the method used
in \cite{CD:Jost:1998} to prove the
classical Noether's theorem. For that we
reparameterize the time (the independent variable $t$) with
a Lipschitzian transformation
\begin{equation*}
[a,b]\ni t\longmapsto \sigma f(\lambda) \in
[\sigma_{a},\sigma_{b}]
\end{equation*}
that satisfies
\begin{equation}
\label{eq:condla} t_{\sigma}^{'} =\frac{dt(\sigma)}{d\sigma}=
f(\lambda) = 1\,\, if\,\, \lambda=0\,.
\end{equation}
In this way one reduces \eqref{Pf} to an autonomous
integral functional
\begin{equation}
\label{eq:tempo} \bar{I}[t(\cdot),q(t(\cdot))] =
\int_{\sigma_{a}}^{\sigma_{b}} \hspace*{-0.2cm}
L\left(t(\sigma),q(t(\sigma)),
{_{\sigma_{a}}^{RC}D_{\sigma_{b}}^{\alpha}q(t(\sigma))}\right)t_{\sigma}^{'}
d\sigma ,
\end{equation}
where $t(\sigma_{a}) = a$ and $t(\sigma_{b}) = b$. Using the
definitions and properties of fractional derivatives given in
Section~\ref{sec:fdRL}, we get successively that
\begin{equation*}
\begin{split}
_{\sigma_{a}}^{RC}D_{\sigma_{b}}^{\alpha}q(t(\sigma))
&= \frac{1}{\Gamma(n-\alpha)}
\int_{\frac{a}{f(\lambda)}}^{\frac{b}{f(\lambda)}}\left|{\sigma
f(\lambda)}-\theta\right|^{n-\alpha-1}\left(\frac{d}{d\theta(\sigma)}\right)^{n}
q\left(\theta f^{-1}(\lambda)\right)d\theta    \\
&= \frac{(t_{\sigma}^{'})^{-\alpha}}{\Gamma(n-\alpha)}
\int_{\frac{a}{(t_{\sigma}^{'})^{2}}}^{\frac{b}{(t_{\sigma}^{'})^{2}}}
|\sigma-s|^{n-\alpha-1}\left(\frac{d}{ds}\right)^{n}
q(s)ds  \\
&= (t_{\sigma}^{'})^{-\alpha}\,\,{_\chi
^{RC}D_{\omega}^{\alpha}}\,q(\sigma) \quad
\left(\chi=\frac{a}{(t_{\sigma}^{'})^{2}}\,, \,\,\,\omega
=\frac{b}{(t_{\sigma}^{'})^{2}}\right)\,.
\end{split}
\end{equation*}
We then have
\begin{equation*}
\begin{split}
\bar{I}[t(\cdot),q(t(\cdot))]
&= \int_{\sigma_{a}}^{\sigma_{b}}
L\left(t(\sigma),q(t(\sigma)),(t_{\sigma}^{'})^{-\alpha}\,\,{_\chi
^{RC}D_{\omega}^{\alpha}}\,q(\sigma)\,q(\sigma)\right) t_{\sigma}^{'} d\sigma \\
&\doteq \int_{\sigma_{a}}^{\sigma_{b}}
\bar{L}_{f}\left(t(\sigma),q(t(\sigma)),t_{\sigma}^{'},\,\,{_\chi
^{RC}D_{\omega}^{\alpha}}\,q(\sigma)\right)d\sigma \\
&= \int_a^b L\left(t,q(t),{_a^{RC}D_b^\alpha} q(t)\right) dt \\
&= I[q(\cdot)] \, .
\end{split}
\end{equation*}
If the integral functional \eqref{Pf} is invariant in the sense of
Definition~\ref{def:invadf}, then the integral functional
\eqref{eq:tempo} is invariant in the sense of
Definition~\ref{def:inv1:MR}. It follows from
Theorem~\ref{theo:tnadf1} that
\begin{equation}
\label{eq:tnadf2}
\mathcal{D}_{t}^{\alpha}\left[\partial_{4}\bar{L}_{f},\xi\right]
+ \mathcal{D}_{t}^{\alpha}\left[\frac{\partial}{\partial t'_\sigma}
\bar{L}_{f},\tau\right] = 0
\end{equation}
is a fractional conserved law in the sense of Riesz-Caputo.
For $\lambda = 0$ the condition \eqref{eq:condla} allows us to
write that
\begin{equation*}
_\chi^{RC}D_\omega^\alpha q(\sigma) = _a^{RC}D_b^\alpha q(t) \, ,
\end{equation*}
and therefore we get that
\begin{equation}
\label{eq:prfMR:q1}
\partial_{4}\bar{L}_{f}=\partial_{3} L
\end{equation}
and
\begin{equation}
\label{eq:prfMR:q2}
\begin{split}
 \frac{\partial}{\partial t'_\sigma} \bar{L}_{f}
&=
\partial_{4}{\bar{L}_{f}} \cdot \frac{\partial}{\partial
t_{\sigma}^{'}}\left[
\frac{(t_{\sigma}^{'})^{-\alpha}}{\Gamma(n-\alpha)}
\int_\chi^{\omega}
|\sigma-s|^{n-\alpha-1}\left(\frac{d}{ds}\right)^{n}q(s)\,ds\right]t_{\sigma}^{'}
 + L \\
&=
\partial_{4}\bar{L}_{f}\cdot\left[\frac{-\alpha(t_{\sigma}^{'})^{-\alpha-1}}
{\Gamma(n-\alpha)} \int_\chi^{\omega}
|\sigma-s|^{n-\alpha-1}\left(\frac{d}{ds}\right)^{n}q(s)\,ds\right]t_{\sigma}^{'}
 + L \\
&= -\alpha\partial_{3} L\cdot{_{a}^{RC}D_{b}^{\alpha}}q
 + L \, .
\end{split}
\end{equation}
Substituting the quantities \eqref{eq:prfMR:q1} and
\eqref{eq:prfMR:q2} into \eqref{eq:tnadf2}, we obtain the
fractional conservation law in the sense of Riesz-Caputo
\eqref{eq:tndf}.
\end{proof}

%--------------------------------------

\subsection{Optimal control of Riesz-Caputo fractional systems}
\label{sub:sec:OC}

We now adopt the Hamiltonian formalism in order to generalize the
Noether type results found in
\cite{CD:Djukic:1972,CD:JMS:Torres:2002a} for the more general context
of fractional optimal control in the sense of Riesz-Caputo.
For this, we make use of our Noether's Theorem~\ref{theo:TNfRL} and
the standard Lagrange multiplier technique (\textrm{cf.}
\cite{CD:Djukic:1972}). The fractional optimal control problem in
the sense of Riesz-Caputo is introduced, without loss of
generality, in Lagrange form:
\begin{equation}
\label{eq:COA} I[q(\cdot),u(\cdot)] = \int_a^b
L\left(t,q(t),u(t)\right) dt \longrightarrow \min \, ,
\end{equation}
subject to the fractional
differential system
\begin{equation}
\label{eq:sitRL}
 _a^{RC}D_a^\alpha
q(t)=\varphi\left(t,q(t),u(t)\right)
\end{equation}
and initial condition
\begin{equation}
\label{eq:COIRL}
q(a)=q_a\, .
\end{equation}
The Lagrangian $L :[a,b] \times \mathbb{R}^{n}\times
\mathbb{R}^{m} \rightarrow \mathbb{R}$ and the fractional velocity vector
$\varphi:[a,b] \times \mathbb{R}^{n}\times \mathbb{R}^m\rightarrow
\mathbb{R}^{n}$ are assumed to be functions of class $C^{1}$ with
respect to all their arguments. We also assume, without loss of
generality, that $0<\alpha\leq1$. In conformity with the calculus
of variations, we are considering that the control functions
$u(\cdot)$ take values on an open set of $\mathbb{R}^m$.

\begin{definition}
The fractional differential system
\eqref{eq:sitRL} is said to be a
\emph{fractional control system
in the sense of Riesz-Caputo}.
\end{definition}

\begin{remark}
In the particular case $\alpha=1$ the problem \eqref{eq:COA}--\eqref{eq:COIRL} is
reduced to the classical optimal control problem
\begin{gather}
I[q(\cdot),u(\cdot)] = \int_a^b
L\left(t,q(t),u(t)\right) dt \longrightarrow \min \, ,\label{eq:ccf}\\
\dot{q}(t)=\varphi\left(t,q(t),u(t)\right)\, ,
\quad q(a)=q_a\, . \label{eq:ccf1}
\end{gather}
\end{remark}

\begin{remark}
\label{rem:cv:pc}
The fractional functional of the calculus of variations in the sense of
Riesz-Caputo \eqref{Pf} is obtained from
\eqref{eq:COA}--\eqref{eq:sitRL} choosing
$\varphi(t,q,u)=u$.
\end{remark}

\begin{definition}[fractional process in the sense of Riesz-Caputo]
An admissible pair $(q(\cdot),u(\cdot))$ which satisfies the
fractional control system \eqref{eq:sitRL} of the
fractional optimal control problem \eqref{eq:COA}--\eqref{eq:COIRL},
$t \in [a,b]$, is said to be a \emph{fractional process in the sense of Riesz-Caputo}.
\end{definition}

\begin{theorem}[\cite{CD:Agrawal:2007}]
\label{th:AG} If $(q(\cdot),u(\cdot))$ is a fractional process
of problem \eqref{eq:COA}--\eqref{eq:COIRL} in the sense of Riesz-Caputo,
then there exists a co-vector function $p(\cdot)\in PC^{1}([a,b];\mathbb{R}^{n})$
such that for all $t\in [a,b]$ the
triple $(q(\cdot),u(\cdot),p(\cdot))$ satisfy the following
conditions:
\begin{itemize}
\item the Hamiltonian system
\begin{equation*}
\label{eq:HamRL}
\begin{cases}
_a^{RC}D_b^\alpha q(t) =\partial_4 {\cal H}(t, q(t), u(t),p(t)) \, , \\
_a^RD_b^\alpha p(t) = -\partial_2{\cal H}(t,q(t),u(t), p(t)) \, ;
\end{cases}
\end{equation*}
\item the stationary condition

\begin{equation*}
 \partial_3 {\cal H}(t, q(t), u(t), p(t))=0 \, ;
\end{equation*}
\end{itemize}
where the Hamiltonian ${\cal H}$ is given by
\begin{equation}
\label{eq:HL}
{\cal H}\left(t,q,u,p\right) = L\left(t,q,u\right)
+ p \cdot \varphi\left(t,q,u\right) \, .
\end{equation}
\end{theorem}

\begin{definition}[fractional Pontryagin extremal in the sense of Riesz-Caputo]
\label{def:extPontRL} A triple $(q(\cdot),u(\cdot),p(\cdot))$
satisfying Theorem~\ref{th:AG} will be called a \emph{fractional
Pontryagin extremal in the sense of Riesz-Caputo}.
\end{definition}

\begin{remark}
In the case of the fractional calculus
of variations in the sense of Riesz-Caputo
one has $\varphi(t,q,u)=u$ (Remark~\ref{rem:cv:pc})
and ${\cal H} = L + p \cdot u$. From Theorem~\ref{th:AG} we get
$_a^{RC}D_a^\alpha q = u$ and
$_a^RD_b^\alpha  p  = -\partial_2 L$ from the
Hamiltonian system, and from
the stationary condition
$\partial_3 {\cal H} = 0$ it follows that
$p= - \partial_3 L$, thus
${_a^RD_b^\alpha}  p= -_a^RD_b^\alpha \partial_3 L$.
Comparing both expressions for $_a^RD_b^\alpha  p$, we arrive to
the fractional Euler-Lagrange equations \eqref{eq:eldf}:
$\partial_2 L = {_a^RD_b^\alpha} \partial_3 L$.
\end{remark}

Minimizing \eqref{eq:COA} subject to
\eqref{eq:sitRL} is equivalent,
by the Lagrange multiplier rule,
to minimize
\begin{equation}
\label{eq:COA1} J[q(\cdot),u(\cdot),p(\cdot)] = \int_a^b
\left[{\cal H}\left(t,q(t),u(t),p(t)\right)-p(t) \cdot
{_a^{RC}D}_a^\alpha q(t)\right]dt
\end{equation}
with ${\cal H}$ given by \eqref{eq:HL}.

\begin{remark}
Theorem~\ref{th:AG} is easily proved applying the
optimality condition \eqref{eq:eldf}
to the equivalent functional \eqref{eq:COA1}.
\end{remark}

The notion of variational invariance for \eqref{eq:COA}--\eqref{eq:sitRL}
is defined with the help of the augmented functional \eqref{eq:COA1}.

\begin{definition}[variational invariance of
\eqref{eq:COA}--\eqref{eq:sitRL}]
\label{def:inv:gt1} We say that the
integral functional \eqref{eq:COA1}
is invariant under the one-parameter
family of infinitesimal transformations
\begin{equation}
\label{eq:trf:inf}
\begin{cases}
\bar{t} = t+\varepsilon\tau(t, q(t), u(t), p(t)) + o(\varepsilon) \, , \\
\bar{q}(t) = q(t)+\varepsilon\xi(t, q(t), u(t), p(t)) + o(\varepsilon) \, , \\
\bar{u}(t) = u(t)+\varepsilon\varrho(t, q(t), u(t), p(t)) + o(\varepsilon) \, , \\
\bar{p}(t) = p(t)+\varepsilon\varsigma(t, q(t), u(t), p(t))+ o(\varepsilon) \, , \\
\end{cases}
\end{equation}
if
\begin{multline}
\label{eq:condInv}
\left[{\cal H}(\bar{t},\bar{q}(\bar{t}),\bar{u}(\bar{t}),\bar{p}(\bar{t}))
-\bar{p}(\bar{t}) \cdot  {_{\bar{a}}^{RC}D_{\bar{b}}}^\alpha
\bar{q}(\bar{t})\right] d\bar{t} \\
=\left[{\cal H}(t,q(t),u(t),p(t))-p(t)
\cdot {_a^{RC}D_b^\alpha} q(t)\right] dt \, .
\end{multline}
\end{definition}

The next theorem provides us with an extension of Noether's
theorem to the wider fractional context of
optimal control in the sense of Riesz-Caputo.

\begin{theorem}[Noether's fractional theorem in Hamiltonian form]
\label{thm:mainResult:FDA06} If \eqref{eq:COA}--\eqref{eq:sitRL}
is variationally invariant, in the sense of Definition~\ref{def:inv:gt1}, then
\begin{multline}
\label{eq:tndf:CO}
\mathcal{D}_{t}^{\alpha}\left[{\cal H}(t,q(t),u(t),p(t)) - \left(1 -
\alpha\right) p(t) \cdot {_a^{RC}D_b^\alpha} q(t),
\tau(t,q(t))\right] \\
- \mathcal{D}_{t}^{\alpha}\left[p(t), \xi(t,q(t))\right] = 0
\end{multline}
along any fractional Pontryagin extremal
$(q(\cdot),u(\cdot),p(\cdot))$ of problem
\eqref{eq:COA}--\eqref{eq:COIRL}.
\end{theorem}
\begin{proof}
The fractional conservation law
\eqref{eq:tndf:CO} in the sense of Riesz-Caputo
is obtained by applying Theorem~\ref{theo:TNfRL}
to the equivalent functional \eqref{eq:COA1}.
\end{proof}

\begin{remark}
In the particular case $\alpha=1$ the fractional optimal control
problem \eqref{eq:COA}--\eqref{eq:COIRL} is reduced to the
standard optimal control problem \eqref{eq:ccf}--\eqref{eq:ccf1}.
In this situation one gets from Theorem~\ref{thm:mainResult:FDA06}
the Noether-type theorem associated with the classical optimal control
problem \cite{CD:Djukic:1972,CD:JMS:Torres:2002a}:
invariance under a one-parameter family of infinitesimal
transformations \eqref{eq:trf:inf} implies that
\begin{equation}
\label{eq:H19} {\cal H}(t,q(t),u(t),p(t))\tau(t,q(t))-p(t)\cdot
\xi(t,q(t)) = constant
\end{equation}
along all the Pontryagin extremals (we obtain the conservation law
\eqref{eq:H19} by choosing $\alpha=1$ in \eqref{eq:tndf:CO}).
\end{remark}

Theorem~\ref{thm:mainResult:FDA06} gives a new and
interesting result for autonomous fractional variational
problems. Let us consider an autonomous fractional
optimal control problem, \textrm{i.e.},
\eqref{eq:COA} and \eqref{eq:sitRL}
with the Lagrangian $L$ and the fractional
velocity vector $\varphi$ not depending explicitly
on the independent variable $t$:
\begin{gather}\label{eq:FOCP:CO}
I[q(\cdot),u(\cdot)] =\int_a^b L\left(q(t),u(t)\right) dt
\longrightarrow \min \, , \\
_a^{RC}D_b^\alpha
q(t)=\varphi\left(q(t),u(t)\right)\label{eq:FOCP:CO10} \, .
\end{gather}

\begin{corollary}
\label{cor:FOCP:CO}
For the autonomous fractional problem \eqref{eq:FOCP:CO}--\eqref{eq:FOCP:CO10}
\begin{equation*}
_a^RD_b^\alpha \left[
{\cal H}(t,q(t),u(t),p(t)) +
\left(\alpha-1\right) p(t) \cdot {_a^{RC}D_b^\alpha} q(t)\right] = 0
\end{equation*}
along any fractional Pontryagin extremal
$(q(\cdot),u(\cdot),p(\cdot))$.
\end{corollary}

\begin{proof}
As the Hamiltonian ${\cal H}$ does not depend explicitly on the
independent variable $t$, we can easily see that
\eqref{eq:FOCP:CO}--\eqref{eq:FOCP:CO10} is invariant
under translation of the time variable:
the condition of invariance \eqref{eq:condInv} is satisfied with
$\bar{t}(t) = t+\varepsilon$,
$\bar{q}(t) = q(t)$,
$\bar{u}(t) = u(t)$,
and $\bar{p}(t) = p(t)$. Indeed,
given that $d\bar{t} = dt$, the invariance
condition \eqref{eq:condInv} is verified if
${_{\bar{a}}^{RC}D_{\bar{b}}^\alpha} \bar{q}(\bar{t}) =
{_a^{RC}D_b^\alpha} q(t)$. This is true because
\begin{equation*}
\begin{split}
{_{\bar{a}}^{RC}D_{\bar{b}}^\alpha}
\bar{q}(\bar{t})
&= \frac{1}{\Gamma(n-\alpha)}
\int_{\bar{a}}^{\bar{b}} |\bar{t}-\theta|^{n
-\alpha-1}\left(\frac{d}{d\theta}\right)^{n}\bar{q}(\theta)d\theta \\
&= \frac{1}{\Gamma(n-\alpha)}
\int_{a + \varepsilon}^{b+\varepsilon} |t + \varepsilon
-\theta|^{n-\alpha-1}\left(\frac{d}{d\theta}\right)^{n}\bar{q}(\theta)d\theta \\
&=\frac{1}{\Gamma(n-\alpha)} \int_{a}^{b}
|t-s|^{n-\alpha-1}\left(\frac{d}{ds}\right)^{n}
\bar{q}(t+\varepsilon)ds\\
&={_a^{RC}D^\alpha_{b}\bar{q}(t+\varepsilon)}={_a^{RC}D^\alpha_{b}\bar{q}(\bar{t})}\\
&={_a^{RC}D^\alpha_{b}{q}(t)}\, .
\end{split}
\end{equation*}
Using the notation in \eqref{eq:trf:inf} we have
$\tau = 1$ and $\xi=\varrho=\varsigma=0$.
From Theorem~\ref{thm:mainResult:FDA06} we arrive to
the intended conclusion.
\end{proof}

Corollary~\ref{cor:FOCP:CO} asserts that unlike
the classical autonomous problem of optimal control,
for \eqref{eq:FOCP:CO}--\eqref{eq:FOCP:CO10}
the fractional Hamiltonian ${\cal H}$ is not conserved.
Instead of $\frac{d}{dt}\left(H\right)=0$ we have
\begin{equation}
\label{eq:ConsHam:alpha}
_a^RD_b^\alpha \left[
{\cal H} + \left(\alpha-1\right) p(t) \cdot {_a^{RC}D_b^\alpha} q(t) \right] = 0 \, ,
\end{equation}
\textrm{i.e.}, fractional conservation
of the Hamiltonian ${\cal H}$ plus a quantity that depends on
the fractional order $\alpha$ of differentiation.
This seems to be explained by violation of the homogeneity
of space-time caused by the fractional
derivatives, $\alpha\neq 1$. In the particular $\alpha=1$
we obtain from \eqref{eq:ConsHam:alpha} the classical result:
the Hamiltonian ${\cal H}$ is preserved along all the Pontryagin
extremals of \eqref{eq:ccf}--\eqref{eq:ccf1}.

%--------------------------------------

\section{Examples}

To illustrate our results, we consider in this section two
examples where the fractional Lagrangian
does not depend explicitly on the
independent variable $t$ (autonomous case).
In both examples we use our Corollary~\ref{cor:FOCP:CO} to
establish the fractional conservation laws.

\begin{example}
Let us consider the following fractional problem
of the calculus of variations:
\begin{equation*}
I[q(\cdot)] = \frac{1}{2}\int_0^1
\left(_{0}^{RC}D_{1}^{\alpha}q(t) \right)^2dt \longrightarrow \min
\, .
\end{equation*}
The Hamiltonian \eqref{eq:HL} takes the form
${\cal H}=-\frac{1}{2}p^2$.
From Corollary~\ref{cor:FOCP:CO} we conclude that
\begin{equation}
\label{eq:lc:ex1}
_0^RD_1^\alpha \left[\frac{p^2(t)}{2}(1-2\alpha)\right] = 0 \, .
\end{equation}
\end{example}

\begin{example}
Consider now the fractional optimal control problem
\begin{gather*}
 I[q(\cdot)] = \frac{1}{2}\int_0^1 \left[q^2(t)+
u^2(t)\right] dt
\longrightarrow \min \, ,\\
_{0}^{RC}D_{1}^{\alpha}q(t)=-q(t)+u(t) \, , \notag
\end{gather*}
under the initial condition $q(0)=1$.
The Hamiltonian ${\cal H}$ defined by \eqref{eq:HL}
takes the following form:
$${\cal H}=\frac{1}{2}\left(q^2+
u^2\right)+p(-q+u).$$
It follows from our Corollary~\ref{cor:FOCP:CO} that
\begin{equation}
\label{eq:cons:Energ:Ex2}
_0^RD_1^\alpha \left[
\frac{1}{2}\left(q^2+u^2\right)
+\alpha p(-q+u) \right] = 0
\end{equation}
along any fractional Pontryagin extremal
$(q(\cdot),u(\cdot),p(\cdot))$ of the problem.
\end{example}

For $\alpha = 1$ the conservation laws
\eqref{eq:lc:ex1} and \eqref{eq:cons:Energ:Ex2} give
the well known results of \emph{conservation of energy}.

%--------------------------------------

\section{Conclusions and Possible Extensions}

The standard approach to solve fractional problems of the calculus of variations
is to make use of the necessary optimality conditions
given by the fractional Euler-Lagrange equations.
These are, in general, nonlinear fractional differential equations,
very hard to be solved. One way to address the problem is to find conservation laws,
\textrm{i.e.}, quantities which are preserved along the Euler-Lagrange
extremals, and that can be used to simplify the problem at hands.
The main questions addressed in this paper are:
(i) What should be a conservation law in the
fractional setting?  (ii) How to determine such fractional conservation laws?
Answer to (i) is given by substituting the classical
operator $\frac{d}{dt}$ by a new fractional operator $\mathcal{D}_{t}^{\gamma}$
(see Definition~\ref{def:oprl} and Remarks~\ref{rem:oprl}, \ref{rem:22}, and \ref{rem:25});
answer to (ii) is given by a fractional Noether's theorem (Theorem~\ref{theo:TNfRL})
that establishes a relation between the existence
of variational symmetries and the existence of fractional conservation laws.

Our results can also bring insight to the issue of nonlinear
fractional differential equations and fractional calculus of variations
from a different perspective.
Indeed, the present results together with
\cite{MyID:183,SilvaTorres06,withLeitmann} suggest
that Leitmann's direct method can also be adapted to cope with
fractional variational problems.
Having in mind that the fractional Euler-Lagrange equations
are in general difficult or even impossible to be solved analytically,
this is an interesting open question,
with some preliminary results already obtained
in \cite{comRic:SI:Leitmann}.

We believe that there is no approach to fractional calculus
better than all the others. All fractional derivatives
have some advantages and disadvantages.
Several fractional variational problems
have been proposed in the literature. This means that
for a given classical Lagrangian
we have at our disposal several different methods to obtain
the fractional Euler-Lagrange equations
and corresponding Hamiltonians.
The fractional dynamics depend on
the fractional derivatives used to construct the Lagrangian,
thus the existence of several options can be used to treat different physical
systems with different specifications and characteristics.
We trust that application of our results to the fractional dynamics
with Riesz-Caputo derivatives will bring new opportunities.

The main aim of our paper was to prove a Noether theorem -- one of
the most beautiful and useful results of the calculus of variations
-- for fractional variational problems via Riesz-Caputo fractional derivatives.
A Noether theorem is valid along the Euler-Lagrange extremals
of the problems. Since previous Euler-Lagrange equations
obtained in the literature via Riesz-Caputo fractional derivatives
(\textrm{cf.} Theorem~\ref{Thm:FractELeq1})
assume both left and right derivatives
to be of the same order $\alpha$,
we have restricted ourselves to this case.
However, in principle it is possible to consider in the
definition of Riesz-Caputo fractional derivative
a different order for the left and right derivative
and prove a corresponding Euler-Lagrange equation.
With such two-order Euler-Lagrange equation one can then
try to generalize our results to this more general setting.
Even more general, one can try to extend
Theorems~\ref{Thm:FractELeq1} and \ref{th:AG}
and Theorems~\ref{theo:TNfRL} and \ref{thm:mainResult:FDA06}
by using generalized Erd\'{e}lyi-Kober fractional integrals and derivatives
\cite{1:ref:2ndrv,2:ref:2ndrv,Virginia94,Virginia08}. This is highly interesting
because of recent applications within the framework of micro structure
\cite{3:ref:2ndrv,4:ref:2ndrv,5:ref:2ndrv,6:ref:2ndrv}.
These are some interesting possibilities for future work.

\bigskip

We are grateful to two anonymous referees for
several good ideas of future research
and encouragement words.

%--------------------------------------

%--------------------------------------

\end{document}